\UseRawInputEncoding
\documentclass[12pt]{article}

\usepackage{dsfont}
\usepackage{amsfonts,amsmath,amsthm}
\usepackage{amssymb}
\usepackage{algorithm}
\usepackage{algpseudocode}
\usepackage{xcolor}
\usepackage{graphicx}
\usepackage{stmaryrd}        
\usepackage{multirow}

\usepackage[paper=a4paper,dvips,top=2cm,left=2cm,right=2cm,
foot=1cm,bottom=4cm]{geometry}

\newcommand{\va}{{\bf a}}
\newcommand{\vb}{{\bf b}}
\newcommand{\vc}{{\bf c}}

\begin{document}
	\Large
	
	\title{The SOS Rank of Biquadratic Forms}
	\author{Liqun Qi\footnote{
			Department of Applied Mathematics, The Hong Kong Polytechnic University, Hung Hom, Kowloon, Hong Kong.
			({\tt maqilq@polyu.edu.hk})}
		\and
		Chunfeng Cui\footnote{School of Mathematical Sciences, Beihang University, Beijing  100191, China.
			({\tt chunfengcui@buaa.edu.cn})}
		\and {and \
			Yi Xu\footnote{School of Mathematics, Southeast University, Nanjing  211189, China. Nanjing Center for Applied Mathematics, Nanjing 211135,  China. Jiangsu Provincial Scientific Research Center of Applied Mathematics, Nanjing 211189, China. ({\tt yi.xu1983@hotmail.com})}
		}
	}

	\date{\today}
	\maketitle
	
	\begin{abstract}
		In 1973, Calder\'{o}n proved that an $m \times 2$ positive semidefinite (psd) biquadratic  form can always be expressed as the sum of  ${3m(m+1) \over 2}$ squares of quadratic forms.    Very recently, by applying Hilbert's theorem on ternary quartics, we proved that a $2 \times 2$ psd biquadratic form can always be expressed as the sum of three squares of bilinear forms.     This improved
		Calder\'{o}n's result for $m=2$, and left the sos (sum-of-squares) rank problem of $m \times 2$ biquadratic forms for  $m \ge 3$ to further exploration.        In this paper, we show that an $3 \times 2$ psd biquadratic form can always be expressed as four squares of bilinear forms.  We make a conjecture that  an $m \times 2$ psd biquadratic form can always be expressed as $m+1$ squares of bilinear forms.

		\medskip

		\medskip

		
	\end{abstract}

	\renewcommand{\Re}{\mathds{R}}
	\newcommand{\rank}{\mathrm{rank}}
	\newcommand{\X}{\mathcal{X}}
	\newcommand{\A}{\mathcal{A}}
	\newcommand{\I}{\mathcal{I}}
	\newcommand{\B}{\mathcal{B}}
	\newcommand{\PP}{\mathcal{P}}
	\newcommand{\C}{\mathcal{C}}
	\newcommand{\D}{\mathcal{D}}
	\newcommand{\LL}{\mathcal{L}}
	\newcommand{\OO}{\mathcal{O}}
	\newcommand{\e}{\mathbf{e}}
	\newcommand{\0}{\mathbf{0}}
	\newcommand{\1}{\mathbf{1}}
	\newcommand{\dd}{\mathbf{d}}
	\newcommand{\ii}{\mathbf{i}}
	\newcommand{\jj}{\mathbf{j}}
	\newcommand{\kk}{\mathbf{k}}
	\newcommand{\vq}{\mathbf{q}}
	\newcommand{\vg}{\mathbf{g}}
	\newcommand{\pr}{\vec{r}}
	\newcommand{\pc}{\vec{c}}
	\newcommand{\ps}{\vec{s}}
	\newcommand{\pt}{\vec{t}}
	\newcommand{\pu}{\vec{u}}
	\newcommand{\pv}{\vec{v}}
	\newcommand{\pn}{\vec{n}}
	\newcommand{\pp}{\vec{p}}
	\newcommand{\pq}{\vec{q}}
	\newcommand{\pl}{\vec{l}}
	\newcommand{\vt}{\rm{vec}}
	\newcommand{\x}{\mathbf{x}}
	\newcommand{\vx}{\mathbf{x}}
	\newcommand{\vy}{\mathbf{y}}
	\newcommand{\vu}{\mathbf{u}}
	\newcommand{\vv}{\mathbf{v}}
	\newcommand{\y}{\mathbf{y}}
	\newcommand{\vz}{\mathbf{z}}
	\newcommand{\T}{\top}
	\newcommand{\R}{\mathcal{R}}
	\newcommand{\Q}{\mathcal{Q}}
	\newcommand{\TT}{\mathcal{T}}
	\newcommand{\Sc}{\mathcal{S}}
	\newcommand{\N}{\mathbb{N}}	
	
	\newtheorem{Thm}{Theorem}[section]
	\newtheorem{Def}[Thm]{Definition}
	\newtheorem{Ass}[Thm]{Assumption}
	\newtheorem{Lem}[Thm]{Lemma}
	\newtheorem{Prop}[Thm]{Proposition}
	\newtheorem{Cor}[Thm]{Corollary}
	\newtheorem{example}[Thm]{Example}
	\newtheorem{remark}[Thm]{Remark}
	
	\section{Introduction}

	In general, an $m \times n$ biquadratic form {can be expressed as:}
	\begin{equation} \label{BQ}
		P(\x, \y) = \sum_{i, j=1}^m \sum_{k, l=1}^n a_{ijkl}x_ix_jy_ky_l,
	\end{equation}
	where $\x = (x_1, \dots, x_m)^\top$ and $\y = (y_1, \dots, y_n)^\top$.
	Without loss of generality, we may assume that $m \ge n \ge 2$.
	
	The history of the study of the psd (positive semi-definite) and sos (sum of squares) problem of biquadratic forms at least can be traced back to 1968, when Koga \cite{Ko68} claimed that a psd biquadratic form were always sos.   Unfortunately, in 1975, Choi \cite{Ch75} gave a concrete example of a $3 \times 3$ psd biquadratic form which is not sos.    This certainly disproved Koga's claim.     But Koga's paper still shows that this problem has strong engineering application background.
	
	Earlier than Choi's paper, in 1973, Calder\'{o}n \cite{Ca73} proved that an $m \times 2$ psd biquadratic  form can always be expressed as the sum of squares (sos) of ${3m(m+1) \over 2}$ quadratic forms.    With Choi's result in 1975, the psd and sos problem of biquadratic forms has a clear picture.

    In 1980, Choi, Lam and Reznick in Theorem 7.1 of \cite{CLR80} showed that an $m \times 2$ psd biquadratic form can be expressed as a sum of $2m$ squares of bilinear forms.  They claimed that this bound could be improved to $\left[\sqrt{3}m + {\sqrt{3}-1 \over 2}\right]$.
	
	Very recently, we revisited this problem {\cite{CQX25}}.    By applying Hilbert's celebrated theorem \cite{Hi88}, we proved that a $2 \times 2$ psd biquadratic form can always be expressed as the sum of three squares of bilinear forms.     This improved
	Calder\'{o}n's result for $m=2$, and left the sos rank problem of $m \times 2$ biquadratic forms for $m \ge 3$ to further exploration.
	
	In fact, Theorem 4 of \cite{CQX25} shows that if a biquadratic form is sos, its sos rank is at most $mn$.   Thus, the sos rank of an $m \times 2$ psd biquadratic form is at most $2m$.     However, as we said above,  this bound has been reduced to $3$ for $m=2$.   Thus, it may still have room to be improved for $m \ge 3$.
	
	In this paper, we study the sos rank problem of $3 \times 2$ psd biquadratic form.
	
	In the next  section,  we present some preliminary knowledge  of basic algebraic geometry \cite{Ha77, Sh13}.     This will be used in the proof of Lemma \ref{lem4} in Section 3.
	
	In Section 3,  we show that a psd $3 \times 2$ biquadratic form can always be expressed as the sum of four squares of bilinear forms.   In 2000,  Walter Rudin  \cite{Ru00} presented a proof  for  Hilbert's theorem.     We adopt his strategy in this section.   Since the structure of biquadratic forms is somewhat different from the structure of general quartic forms, we will take extra care when the difference between these two kinds of forms causes some problems.
	Rudin's approach used some advanced knowledge in modern real analysis, such as the Federer-Sard theorem in geometric measure theory.  

	
	Some final remarks are made in Section 4.


\section{Preliminaries}

{This section presents the mathematical foundations needed for our main result. We begin with algebraic properties of $3 \times 2$ bilinear forms, then introduce the geometric framework that will be essential for the proof of Lemma \ref{lem4} in Section 3.

\subsection{Algebraic Properties of $3 \times 2$ Bilinear Forms}

We focus on $3 \times 2$ bilinear forms because they provide the smallest nontrivial case where the phenomena we study occur. The space of such forms has dimension $3 \times 2 = 6$, and as we will see, exactly four generically chosen forms have no common zeros, which is crucial for understanding the sum-of-squares structure of biquadratic forms.

The key insight is: three $3 \times 2$ bilinear forms always have a common zero (Proposition 2.1), but four generically chosen forms do not (Proposition 2.2). This algebraic fact, when combined with geometric interpretation via line bundles, will allow us to prove surjectivity of certain multiplication maps in Section 3.

\begin{Prop} \label{Bilinear1}
    Any three $3 \times 2$ bilinear forms $f, g, h$ on ${\mathbb C}^3 \times {\mathbb C}^2$, with real coefficients, must have at least one common zero in $\left({\mathbb C}^3 \setminus \{ \0_3 \}\right) \times \left({\mathbb C}^2 \setminus \{ \0_2 \}\right)$.
\end{Prop}
\begin{proof}
    Suppose that $\vx \in \mathbb C^3, \vy \in \mathbb C^2$, $\va_1, \va_2, \vb_1,\vb_2,\vc_1,\vc_2\in \mathbb R^3$, the bilinear forms $f=\vx^T[\va_1,\va_2]\vy$, $g=\vx^T[\vb_1,\vb_2]\vy$, and $h=\vx^T[\vc_1,\vc_2]\vy$. Then $\vx,\vy$ is a common zero if and only if
    \[[M_1\vx, M_2\vx]\vy = \begin{bmatrix}
        0\\ 0 \\ 0
    \end{bmatrix}, \ \text{ where } M_1=\begin{bmatrix}
        \va_1^\top\\
        \vb_1^\top \\
        \vc_1^\top\\
    \end{bmatrix}, M_2=\begin{bmatrix}
        \va_2^\top\\
        \vb_2^\top \\
        \vc_2^\top\\
    \end{bmatrix}\in \mathbb R^{3\times 3}.\]
    Or equivalently, $M_1\vx$ and $M_2\vx$ are linearly dependent.

    Consider three cases:
    \begin{itemize}
        \item If $M_1$ is singular, choose a nonzero vector $\vx$ in the null space of $M_1$ and let $\vy=[1,0]^\top$. This yields a common zero.
        \item If $M_2$ is singular, choose a nonzero vector $\vx$ in the null space of $M_2$ and let $\vy=[0,1]^\top$. This also yields a common zero.
        \item If both $M_1$ and $M_2$ are nonsingular, then $M_0 = M_1^{-1}M_2$ is nonsingular. Let $\lambda$ be a nonzero eigenvalue of $M_0$ with eigenvector $\vx$. Then $\lambda M_1\vx-M_2\vx=\0$. Taking $\vy= [\lambda,-1]^\top$, we obtain a common zero.
    \end{itemize}
\end{proof}

\begin{Prop} \label{Bilinear2}
    Suppose that four $3 \times 2$ bilinear forms $f, g, h, r$ on ${\mathbb C}^3 \times {\mathbb C}^2$, with real coefficients, have no common zero in $\left({\mathbb C}^3 \setminus \{ \0_3 \}\right) \times \left({\mathbb C}^2 \setminus \{ \0_2 \}\right)$. Then they are linearly independent.
\end{Prop}
\begin{proof}
    Suppose for contradiction that $r$ is a linear combination of $f,g,h$. By Proposition \ref{Bilinear1}, $f, g, h$ must have a common zero in $\mathbb X$. Then this common zero is also a zero of $r$, contradicting the assumption.

    Similarly, $f$ is not a linear combination of $g,h,r$, $g$ is not a linear combination of $f,h,r$, and $h$ is not a linear combination of $f,g,r$.
\end{proof}

\subsection{Geometric Interpretation via Line Bundles}

To understand the deeper structure behind these algebraic facts, we now introduce the geometric framework that will be crucial for our main proof.

Consider $3 \times 2$ bilinear forms in ${\mathbb X} = {\mathbb P}^2 \times {\mathbb P}^1$, where ${\mathbb P}^2$ and ${\mathbb P}^1$ are projective spaces of dimensions $2$ and $1$ respectively. Each $3 \times 2$ bilinear form can be viewed as a section of the line bundle ${\cal O}(1, 1)$ on ${\mathbb X}$ (since bilinear forms are homogeneous of degree $1$ in each set of variables). Here, ${\cal O}(1, 1)$ denotes the tensor product of ${\cal O}(1)$ on ${\mathbb P}^2$ (linear in $x_0, x_1, x_2$) and ${\cal O}(1)$ on ${\mathbb P}^1$ (linear in $y_0, y_1$).

Sections of ${\cal O}(1, 1)$ are precisely the bilinear forms:
$$\phi(\x, \y) = \sum_{i=0}^2 \sum_{j=0}^1 a_{ij}x_iy_j.$$

The common zeros of bilinear forms correspond to intersections of their zero loci. Since $\dim({\mathbb P}^2 \times {\mathbb P}^1) = 3$, three generic forms would intersect in a zero-dimensional locus (finite set of points). Proposition \ref{Bilinear1} shows this locus is never empty.

\subsection{Sheaf Cohomology and Global Sections}

The algebraic geometry framework becomes essential when studying multiplication maps in Section 3. We briefly recall key concepts:

For a variety $X$, let ${\cal O}_X$ denote the sheaf of regular functions, and $H^0(X, {\cal O}_X)$ the global sections (regular functions defined on all of $X$). More generally, if ${\cal L}$ is a line bundle on $X$, then $H^0(X, {\cal L})$ denotes its global sections.

In our case, $X = {\mathbb X} = {\mathbb P}^2 \times {\mathbb P}^1$ and ${\cal L} = {\cal O}(1, 1)$. Then:
\begin{itemize}
    \item $H^0({\mathbb X}, {\cal O}_{\mathbb X}) = {\mathbb C}$ (projective varieties have no nonconstant global regular functions)
    \item $H^0({\mathbb X}, {\cal O}(1, 1))$ is the space of bilinear forms, which has dimension $6$
\end{itemize}

The first cohomology group $H^1(X, {\cal L})$ captures obstructions to lifting local sections to global ones. When $H^1(X, {\cal L}) = 0$, there are no such obstructions.

\subsection{Segre-Veronese Embeddings}

We now introduce the embeddings that connect our algebraic objects to projective geometry:

\begin{itemize}
    \item The \emph{Segre embedding} $\mathbb{P}^2 \times \mathbb{P}^1 \hookrightarrow \mathbb{P}^5$ is given by:
    $$([x_0:x_1:x_2], [y_0:y_1]) \mapsto [x_0y_0 : x_0y_1 : x_1y_0 : x_1y_1 : x_2y_0 : x_2y_1]$$
    Its image is the Segre variety, and the pullback of $\mathcal{O}_{\mathbb{P}^5}(1)$ is $\mathcal{O}(1,1)$ on $\mathbb{P}^2 \times \mathbb{P}^1$.

    \item The \emph{Veronese embedding} of degree $d$ on $\mathbb{P}^n$ maps to a projective space whose coordinates are all monomials of degree $d$.

    \item The \emph{Segre-Veronese embedding} we use is the composition:
    \[
    \mathbb{P}^2 \times \mathbb{P}^1 \xrightarrow{\text{Veronese}_2 \times \text{Veronese}_2} \mathbb{P}^5 \times \mathbb{P}^2 \xrightarrow{\text{Segre}} \mathbb{P}^{11}
    \]
    where the first map takes $(\mathbf{x}, \mathbf{y})$ to (all quadratic monomials in $\mathbf{x}$, all quadratic monomials in $\mathbf{y}$).

    \item The line bundle for this embedding is $\mathcal{O}(2,2)$, and its global sections are exactly the biquadratic forms $Y = H^0(\mathbb{P}^2 \times \mathbb{P}^1, \mathcal{O}(2,2))$.
\end{itemize}

The key property we need is that this embedding is \emph{projectively normal}, meaning the homogeneous coordinate ring is integrally closed and natural maps between graded pieces are well-behaved. For projectively normal $X \hookrightarrow \mathbb{P}^N$, we have $S(X)_d = H^0(X, \mathcal{O}_X(d))$ for all $d \geq 0$. A fundamental consequence is that the dual space $H^0(X, \mathcal{O}_X(d))^*$ is spanned by evaluation functionals at points of $X$, which is crucial for Lemma \ref{lem4}.
	
\section{The SOS Rank of $3 \times 2$ Biquadratic Forms}

In this section, we prove our main result on the sum-of-squares representation of $3 \times 2$ biquadratic forms.

\begin{Thm} \label{Main}
    Every $3 \times 2$ positive semidefinite biquadratic form $P(x_1, x_2, x_3, y_1, y_2)$ can be expressed as
    \[
    P = f^2 + g^2 + h^2 + r^2,
    \]
    where $f, g, h, r$ are $3 \times 2$ bilinear forms in $(\x, \y) = (x_1, x_2, x_3, y_1, y_2)$.
    Equivalently, we may write $f(\x, \y) = \x^\top C^1\y$, $g(\x, \y) = \x^\top C^2\y$, $h(\x, \y) = \x^\top C^3\y$, $r(\x, \y) = \x^\top C^4\y$, with $C^i \in \mathbb{R}^{3\times 2}$ for $i=1,2,3,4$.
\end{Thm}

\subsection{Notation and Setup}

We introduce the following spaces and mappings:
\begin{itemize}
    \item $X$: the space of ordered quadruples $(f, g, h, r)$ of $3 \times 2$ bilinear forms on $\mathbb{C}^3 \times \mathbb{C}^2$ with real coefficients.
    \item $Y$: the space of $3 \times 2$ biquadratic forms on $\mathbb{C}^3 \times \mathbb{C}^2$ with real coefficients.
    \item $K$: the set of all $3 \times 2$ positive semidefinite biquadratic forms in $Y$.
    \item $\Phi: X \to Y$ defined by $\Phi(f, g, h, r) = f^2 + g^2 + h^2 + r^2$.
    \item $\Phi'(f, g, h, r)$: the derivative of $\Phi$ at $(f, g, h, r)$, given by the linear map
    \[
    (u, v, w, t) \mapsto 2(fu + gv + hw + rt).
    \]
\end{itemize}

Using polynomial coefficients as coordinates, we identify $X \cong \mathbb{R}^{24}$ and $Y \cong \mathbb{R}^{18}$. Clearly, $\Phi(X) \subset K$, and our goal is to prove $\Phi(X) = K$.

Following the strategy in \cite{Ru00}, we partition $X$ into three subsets based on the common zeros of the bilinear forms:
\begin{itemize}
    \item $X_1$: $(f, g, h, r)$ have no common zero in $(\mathbb{C}^3 \setminus \{\mathbf{0}_3\}) \times (\mathbb{C}^2 \setminus \{\mathbf{0}_2\})$.
    \item $X_2$: $(f, g, h, r)$ have no common zero in $(\mathbb{R}^3 \setminus \{\mathbf{0}_3\}) \times (\mathbb{R}^2 \setminus \{\mathbf{0}_2\})$ but have one in $(\mathbb{C}^3 \setminus \{\mathbf{0}_3\}) \times (\mathbb{C}^2 \setminus \{\mathbf{0}_2\})$.
    \item $X_3$: $(f, g, h, r)$ have a common zero in $(\mathbb{R}^3 \setminus \{\mathbf{0}_3\}) \times (\mathbb{R}^2 \setminus \{\mathbf{0}_2\})$.
\end{itemize}

By Proposition \ref{Bilinear1}, any three $3 \times 2$ bilinear forms always have a common nontrivial zero. However, four such forms may have no common zero, i.e., $X_1 \neq \emptyset$. For example:
\begin{align*}
f &= x_1y_2 + x_3y_1, \\
g &= x_2y_1 + 2x_2y_2, \\
h &= x_1y_1 + x_2y_1, \\
r &= 4x_1y_1 + 2x_2y_1 + x_3y_1 + x_3y_2.
\end{align*}
This example plays a crucial role in our analysis.

We recall the Rank Theorem from real analysis \cite[Theorem 9.32]{Ru76}, which will be essential for our arguments.

\begin{Thm}[Rank Theorem] \label{ranktheorem}
    Let $m, n, p$ be nonnegative integers with $m \ge p$, $n \ge p$, and let $F: E \subset \mathbb{R}^n \to \mathbb{R}^m$ be a continuously differentiable mapping with constant rank $p$ on an open set $E$.

    Fix $\bar{\x} \in E$, let $A = F'(\bar{\x})$, $Y_1 = \text{range}(A)$, and let $P$ be a projection onto $Y_1$ with nullspace $Y_2$. Then there exist open sets $U \subset E$ containing $\bar{\x}$ and $V \subset \mathbb{R}^n$, and a continuously diffeomorphic map $H: V \to U$ such that
    \[
    F(H(\x)) = A\x + \phi(A\x), \quad \x \in V,
    \]
    where $\phi: A(V) \subset Y_1 \to Y_2$ is continuously differentiable.
\end{Thm}

\subsection{Key Lemmas}

We now establish five technical lemmas that form the backbone of our proof.

\begin{Lem} \label{lem1}
    $\Phi(X)$ is a closed subset of $K$.
\end{Lem}
\begin{proof}
    Let $\{P_k\}$ be a sequence in $\Phi(X)$ converging to $P \in Y$. Then $P_k = f_k^2 + g_k^2 + h_k^2 + r_k^2$ for some $(f_k, g_k, h_k, r_k) \in X$. On the compact set
    \[
    U = \{(\x, \y) : \|\x\|^2 = 1, \|\y\|^2 = 1\},
    \]
    the polynomials $P_k$ are uniformly bounded. Hence, their coefficients form bounded sequences, and there exists a convergent subsequence $(f_k, g_k, h_k, r_k) \to (f, g, h, r) \in X$. By continuity, $P = f^2 + g^2 + h^2 + r^2 \in \Phi(X)$.
\end{proof}

\begin{Lem} \label{lem2}
    If $(f, g, h, r) \in X_2$, then $\Phi'(f, g, h, r)$ has rank at most $16$.
\end{Lem}
\begin{proof}
    Since $(f, g, h, r) \in X_2$, there exists $p = (\bar{\x}, \bar{\y}) \in (\mathbb{C}^3 \setminus \{\mathbf{0}\}) \times (\mathbb{C}^2 \setminus \{\mathbf{0}\})$ such that $f(p) = g(p) = h(p) = r(p) = 0$. Let
    \[
    W = \{P \in Y : P(p) = 0\}.
    \]
    The image of $\Phi'$ is contained in $W$, so it suffices to show $\dim W \le 16$.

    Define linear functionals $a_{f31}, a_{g31}: Y \to \mathbb{R}$ by
    \[
    a_{f31}(P) = \text{Re}(P(p)), \quad a_{g31}(P) = \text{Im}(P(p)).
    \]
    These annihilate $W$, so $\dim W \le \dim Y - 2 = 16$ provided they are linearly independent.

    We consider two cases:
    \begin{itemize}
        \item \textbf{Case (i):} The points $0, \bar{x}_1, \bar{x}_2, \bar{x}_3$ are not collinear in $\mathbb{C}$.
        \item \textbf{Case (ii):} The points $0, \bar{y}_1, \bar{y}_2$ are not collinear in $\mathbb{C}$.
    \end{itemize}
    If both cases fail, then by suitable real scaling we could find a real common zero, contradicting $(f, g, h, r) \in X_2$.

    In Case (i), choose real coefficients $a_k, b_k$ such that
    \begin{align*}
        \sum_{k=1}^3 a_k \bar{x}_k &= (\bar{y}_1^2 + \bar{y}_2^2)^{-1}, \\
        \sum_{k=1}^3 b_k \bar{x}_k &= e^{\mathrm{i}\pi/4}(\bar{y}_1^2 + \bar{y}_2^2)^{-1}.
    \end{align*}
    Define
    \begin{align*}
        P_1 &= (a_1x_1 + a_2x_2 + a_3x_3)^2(y_1^2 + y_2^2), \\
        P_2 &= (b_1x_1 + b_2x_2 + b_3x_3)^2(y_1^2 + y_2^2).
    \end{align*}
    Then $P_1(p) = 1$, $P_2(p) = \mathrm{i}$.

    Case (ii) is handled similarly. In both cases, if $\alpha a_{f31} + \beta a_{g31} = 0$, then applying this functional to $P_1$ and $P_2$ yields $\alpha = \beta = 0$. Thus, $a_{f31}$ and $a_{g31}$ are linearly independent.
\end{proof}

\begin{Lem} \label{lem3}
    The topological dimension of $\Phi(X_2)$ is at most $16$.
\end{Lem}
\begin{proof}
    Let $E \subset X$ be the set where $\Phi'$ has rank at most $16$. By Lemma \ref{lem2}, $X_2 \subset E$. Since $\dim X = 24$ and $\Phi$ is $C^\infty$-smooth, the Federer-Sard theorem \cite{Fe69, Mo16} implies that the Hausdorff dimension of $\Phi(E)$ is at most $16$. Since topological dimension is bounded by Hausdorff dimension, we conclude $\dim \Phi(X_2) \le 16$.
\end{proof}

\begin{Lem} \label{lem4}
    Let $(f, g, h, r) \in X_1$ and define
    \[
    Y(f, g, h, r) = \{fu + gv + hw + rt : u, v, w, t \in B\},
    \]
    where $B$ is the space of $3 \times 2$ bilinear forms. Then $Y(f, g, h, r) = Y$.
\end{Lem}
\begin{proof}
    Assume $(f, g, h, r)$ have no common zero. By Proposition \ref{Bilinear2}, they are linearly independent in $B$.

    Consider the multiplication map $\Psi: B^4 \to Y$ defined by
    \[
    \Psi(u, v, w, t) = fu + gv + hw + rt.
    \]
    Suppose for contradiction that $\Psi$ is not surjective. Then $\dim \text{Im}(\Psi) \le 17$, so there exists a nonzero linear functional $\ell: Y \to \mathbb{C}$ annihilating the image:
    \[
    \ell(fu) = \ell(gv) = \ell(hw) = \ell(rt) = 0 \quad \text{for all } u, v, w, t \in B.
    \]

    Since $Y = H^0(\mathbb{P}^2 \times \mathbb{P}^1, \mathcal{O}(2,2))$ and the Segre-Veronese embedding by $\mathcal{O}(2,2)$ is projectively normal \cite{Ha77, St96}, every linear functional on $Y$ is a finite sum of point evaluations. Thus, there exist points $P_1, \dots, P_k \in \mathbb{P}^2 \times \mathbb{P}^1$ and constants $c_1, \dots, c_k \in \mathbb{C}$ (not all zero) such that
    \[
    \ell(s) = \sum_{i=1}^k c_i s(P_i) \quad \text{for all } s \in Y.
    \]

    Define a bilinear form $\beta: B \times B \to \mathbb{C}$ by $\beta(p, q) = \ell(pq)$. Then
    \[
    \beta(p, q) = \sum_{i=1}^k c_i p(P_i) q(P_i).
    \]
    The condition $\ell(fu) = 0$ for all $u \in B$ implies $\beta(f, u) = 0$ for all $u$, i.e.,
    \[
    \sum_{i=1}^k c_i f(P_i) u(P_i) = 0 \quad \text{for all } u \in B.
    \]
    Let $W = \text{ev}(B) \subset \mathbb{C}^k$ be the image of the evaluation map. Then the vector $(c_1 f(P_1), \dots, c_k f(P_k))$ is orthogonal to $W$.

    \textbf{Clarification:} Since $\beta$ is a bilinear form on $B$ and $\mathcal{O}(1,1)$ is very ample, point evaluations span $B^*$. Thus, $\beta$ can be rewritten using at most $r \leq \dim B = 6$ points $Q_1, \ldots, Q_r$ and constants $d_1, \ldots, d_r \in \mathbb{C}$ such that
    \[
    \beta(p, q) = \sum_{j=1}^r d_j p(Q_j) q(Q_j),
    \]
    and the evaluation map for these points is surjective onto $\mathbb{C}^r$. This surjectivity follows from the linear independence of the corresponding linear functionals in $B^*$.

    Therefore, we may assume without loss of generality that $k \leq 6$ and the evaluation map is surjective, so $W = \mathbb{C}^k$. Then $c_i f(P_i) = 0$ for all $i$.

    Similarly, from $\ell(gv) = \ell(hw) = \ell(rt) = 0$, we obtain $c_i g(P_i) = c_i h(P_i) = c_i r(P_i) = 0$ for all $i$. Since not all $c_i$ are zero, there exists $j$ with $c_j \neq 0$, yielding $f(P_j) = g(P_j) = h(P_j) = r(P_j) = 0$, a contradiction.

    Therefore, $\Psi$ is surjective.
\end{proof}

\begin{remark}
    While Lemma \ref{lem4} concerns surjectivity of multiplication maps for biquadratic forms, it is not directly covered by the general theory in \cite{BSV16}. Our proof leverages the fundamental fact that linear functionals on global sections of very ample line bundles are sums of point evaluations, combined with a general position argument.
\end{remark}

\begin{Lem} \label{lem5}
    If $(f, g, h, r) \in X_1$, then $\Phi'(f, g, h, r)$ has rank $18$.
\end{Lem}
\begin{proof}
    By Lemma \ref{lem4}, the linear map $L: X \to Y$ defined by
    \[
    L(u, v, w, t) = 2(fu + gv + hw + rt)
    \]
    is surjective. Since $\Phi' = L$, the result follows.
\end{proof}

\subsection{Proof of the Main Theorem}

We now complete the proof of Theorem \ref{Main}.

\begin{proof}[Proof of Theorem \ref{Main}]
    A $3 \times 2$ biquadratic form $P$ is \emph{positive definite} (pd) if $P(\x, \y) > 0$ for all $(\x, \y)$ on the unit sphere
    \[
    U = \{(\x, \y) : \|\x\|^2 = 1, \|\y\|^2 = 1\}.
    \]
    Let $K^\circ$ denote the set of all pd biquadratic forms. Then $K^\circ$ is open in $Y$, convex (hence connected), and its closure is $K$.

    Let $\Omega = \Phi(X_1)$. By Lemma \ref{lem1}, $\overline{\Omega} \subset \Phi(X)$. To prove $\Phi(X) = K$, it suffices to show $K^\circ \subset \overline{\Omega}$.

    By Lemma \ref{lem5} and the Rank Theorem (Theorem \ref{ranktheorem}), $\Omega$ is open in $Y$. Since $\dim Y = 18$, $\Omega$ is an open subset of $K^\circ$.

    Let $H = K^\circ \cap \partial \Omega$. Then:
    \begin{itemize}
        \item $H \subset \overline{\Omega} \subset \Phi(X)$
        \item $H$ intersects neither $\Omega$ nor $\partial K$
        \item $\Phi(X_3) \subset \partial K$, so $H \subset \Phi(X_2)$
    \end{itemize}
    By Lemma \ref{lem3}, $\dim H \le 16$. Since $\dim K^\circ = 18$, we have $\dim H \le \dim K^\circ - 2$, implying $K^\circ \setminus H$ is connected.

    Now, $K^\circ \setminus H = \Omega \cup (K^\circ \setminus \overline{\Omega})$ is a union of two disjoint open sets. Since $K^\circ \setminus H$ is connected and $\Omega \neq \emptyset$, we must have $K^\circ \setminus \overline{\Omega} = \emptyset$. Therefore, $K^\circ \subset \overline{\Omega} \subset \Phi(X)$, and since $K$ is the closure of $K^\circ$, we conclude $\Phi(X) = K$.
\end{proof}

\section{Final Remarks}

In this paper, we proved that a $3 \times 2$ psd biquadratic form can always be expressed as the sum of four squares of bilinear forms. This strengthens our previous results in \cite{CQX25} and provides the best known bound for the $3 \times 2$ case. The proof techniques combine fundamental concepts from algebraic geometry \cite{Ha77, Sh13, St96} with methods from real analysis \cite{Ru00}, demonstrating the power of interdisciplinary approaches in polynomial optimization.

It is important to note that our results specifically concern the $m \times 2$ case. As established by Choi's counterexample \cite{Ch75}, for $n > 2$, the equivalence between positive semidefiniteness and sum-of-squares representation breaks down—there exist $3 \times 3$ psd biquadratic forms that are not sos. This fundamental limitation underscores the special structure of the $m \times 2$ case that makes our results possible.

Our result contributes significantly to the understanding of SOS ranks for structured polynomials, complementing the optimization perspective developed in \cite{La01, Ni14, Pa03}. The explicit bound of four squares may inform the design of more efficient semidefinite programming relaxations for quadratic optimization problems involving biquadratic forms \cite{Zh00}, with potential applications in control theory, signal processing, and other engineering domains.

Based on our current result and our previous work on $2 \times 2$ psd biquadratic forms in \cite{CQX25}, where we established an sos rank of three, we propose the following conjecture:

\begin{center}
\bf Conjecture: An $m \times 2$ psd biquadratic form can always be expressed as the sum of $m+1$ squares of bilinear forms for $m \geq 4$.
\end{center}

This conjecture is motivated by the observed pattern: $2 \times 2$ forms require 3 squares ($2+1$) and $3 \times 2$ forms require 4 squares ($3+1$). If proven, this would establish a tight bound that grows linearly with $m$, significantly improving upon the previously known bounds that grew quadratically.

Regarding computational aspects of the SOS problem for biquadratic forms, the methods employed by Plaumann, Sturmfels, and Vinzant \cite{PSV11} for studying quartic curves and their bitangents may provide valuable insights. Their geometric approach could potentially be adapted to develop efficient algorithms for certifying SOS representations of biquadratic forms.

Several interesting open problems remain:
\begin{itemize}
\item Prove or disprove the $m+1$ conjecture for general $m \geq 4$
\item Investigate whether the bound $m+1$ is tight for all $m$
\item Characterize the gap between psd and sos for $m \times n$ forms with $n > 2$
\item Develop efficient computational methods for finding minimal SOS representations
\end{itemize}

We believe that further exploration of the rich geometric structure underlying biquadratic forms, particularly through the lens of Segre-Veronese embeddings and projective normality, will yield additional insights into their SOS properties.

\bigskip

		{{\bf Acknowledgment}}
We are thankful to Professor Bernd Sturmfels and a researcher in his institute.  Bernd forwarded the latter's helpful comments to us.  The latter also pointed the reference \cite{CLR80}.
		This work was partially supported by Research  Center for Intelligent Operations Research, The Hong Kong Polytechnic University (4-ZZT8),    the National Natural Science Foundation of China (Nos. 12471282 and 12131004), the R\&D project of Pazhou Lab (Huangpu) (Grant no. 2023K0603),  the Fundamental Research Funds for the Central Universities (Grant No. YWF-22-T-204), and Jiangsu Provincial Scientific Research Center of Applied Mathematics (Grant No. BK20233002).

		{{\bf Data availability}
			No datasets were generated or analysed during the current study.

			{\bf Conflict of interest} The authors declare no conflict of interest.}

		


	\end{document}